\documentclass[reqno]{siamltex}

\usepackage{amsmath}
\usepackage{amscd}
\usepackage{amssymb}
\usepackage[usenames]{color}

\usepackage{amsfonts,amsmath,bm}
\usepackage{amssymb}
\usepackage{euscript,graphicx,subfigure,times,version}
\usepackage{fourier}
\newcommand{\bdm}{\begin{displaymath}}
\newcommand{\edm}{\end{displaymath}}
\newcommand{\define}{:=}
\newcommand{\K}{{\mathcal K}}

\newcommand{\dd}[1]{\, {\mathrm d}{#1}}

\newcommand{\difh}[3]{{\frac{\mathrm{d}^{#1}{#2}}{\mathrm{d}\,{#3}^{#1}}}}

\renewcommand{\(}{\left(}
\renewcommand{\)}{\right)}
\newcommand{\N}{{\mathbb N}}

\newcommand{\Rf}{{\mathbb R}}
\newcommand{\C}{{\mathbb C}}
\newcommand{\im}{{\mathrm i}}
\newcommand{\e}{{\mathrm e}}
\newcommand{\half}{\frac{1}{2}}

\excludeversion{hide}

\newtheorem{thm}{Theorem}
\newtheorem{prop}{Proposition}

\newtheorem{rmk}{Remark}

\begin{hide}
\theoremstyle{plain}

\theoremstyle{definition}

\newenvironment{proof}{\noindent{\it Proof}\rm.}{\hfill $\Box$}

\end{hide}

\numberwithin{equation}{section}

\title{On Rayleigh-Type Formulas for a Non-local Boundary Value Problem
Associated with an Integral Operator Commuting with the Laplacian
\thanks{
Received by the editors September 26, 2013; revised XXXX XX, XXXX.; accepted for publication (in revised form) XXXX XX, XXXX.; published electronically XXXX XX, XXXX.
}}

\author{
Lotfi Hermi
\thanks{Department of Mathematics,
University of Arizona,
617 Santa Rita, Tucson, AZ 85721 USA
({\tt hermi@math.arizona.edu}).}
\and 
Naoki Saito
\thanks{Department of Mathematics,
University of California, 
One Shields Avenue,
Davis, CA 95616 USA
({\tt saito@math.ucdavis.edu}).}
}

\begin{document}
\maketitle

\maketitle

\begin{abstract}
In this article we prove the existence, uniqueness, and simplicity of a negative eigenvalue for a class of integral operators whose kernel is of the form $|x-y|^\rho$, $0 < \rho \leq 1$, $x, y \in [-a, a]$. We also provide two different ways of producing recursive formulas for the Rayleigh functions (i.e., recursion formulas for power sums) of the eigenvalues of this integral operator when $\rho=1$, providing means of approximating this negative eigenvalue. These methods offer recursive procedures for dealing with the eigenvalues of a one-dimensional Laplacian with non-local boundary conditions which commutes with an integral operator having a harmonic kernel. The problem emerged in recent work by one of the authors \cite{SAITO-LAPEIG-ACHA}. We also discuss extensions in higher dimensions and links with distance matrices. 
\end{abstract}

\medskip



\medskip




\pagestyle{myheadings} \thispagestyle{plain}
\markboth{L. HERMI AND N. SAITO}{Rayleigh Functions}

\section{Introduction}
\label{sec:intro}

There has been renewed interest, motivated by applications in statistics, machine learning, and mathematical physics, in the spectral properties of integral operators \cite{BapatKirklandNeumann,Bekers,Bobo1,Bobo2,ChoquardStubbe,diaconis,Karoui,Jiang,SAITO-LAPEIG-ACHA}. These operators are usually defined in terms of symmetric distance-like kernels where the focus has recently shifted to questions about spectral embedding, and on establishing connections between empirical operators and their continuous counterparts \cite{Rosasco-Belkin-DeVito}, specifically in the context of manifold learning, with recent activities \cite{Bavaud,Bobo1,Bobo2,diaconis} reviving the theories developed by Schoenberg in the 1930s \cite{Schoenberg1,Schoenberg2,Schoenberg3}, or borrowing techniques from the discrete setting to approximate eigenvalues and eigenfunctions for the continuous counterpart \cite{Bekers,diaconis,Reade}. As a prototype of such integral operators, we consider 
\begin{equation} \label{operator}
\K_{\rho, a} f(x) \define  C_\rho \displaystyle{ \int_{-a}^{a} |x-y|^{\rho} \, f(y) \dd{y}}
\end{equation}
where $a>0$ and
\begin{equation} \label{ct} 
C_\rho \define \frac{\Gamma(-\rho)}{\Gamma\(\frac{1-\rho}{2}\) \, \Gamma\(\frac{1+\rho}{2}\)}= \frac{-1}{2 \Gamma(1+\rho) \sin \frac{\pi \rho}{2}} \quad < 0
\end{equation}
for $0<\rho \le 1$, $C_1=\lim_{\rho \to 1} C_\rho =-1/2$.

The constant $C_\rho$ is motivated by the decomposition of $|x-y|^{\rho}$, due to P\'olya-Szeg\H{o} \cite{PS}, who proved that for $-1 \leq x, y \leq 1$, $-1 < \rho < 1$, with $x \neq y$, $\rho \neq 0$, 
\begin{equation} \label{PS-decomposition}
|x-y|^{\rho} = \frac{\Gamma\(\frac{1+ \rho}{2}\) \, \Gamma\(1-\frac{\rho}{2}\)}{\Gamma\(\frac{1}{2}\)} \, \sum_{n=0}^{\infty} \(1-\frac{2n}{\rho} \) \, P_n^{\(-\frac{\rho}{2}\)} (x) P_n^{\(-\frac{\rho}{2}\)} (y)
\end{equation}
(see Eq.(14) of \cite{PS}, and the comments on p. 29 just before Eq.(18), beginning ``Die Entwicklung (14) \dots'').
They also established the identity
\begin{equation} \label{identity}
\int_{-1}^1 \left(1-x^2\right)^{-\frac{1+\rho}{2}} \, |x-y|^{\rho} P_n^{\(-\frac{\rho}{2}\)} (x) \dd{x} = \frac{\Gamma\(\frac{1-\rho}{2}\) \, \Gamma\(\frac{1+\rho}{2}\)}
{\Gamma(-\rho)} \, \frac{\Gamma(n-\rho)}{\Gamma(n+1)} \, P_n^{\(-\frac{\rho}{2}\)} (y).
\end{equation}
Here $P_n^{(\nu)}(x)$ denotes the ultraspherical (or Gegenbauer) polynomials.
In this article we use the classical notation for Gegenbauer polynomials 
rather than the more modern $C_n^{(\nu)}(x)$ found in 
e.g., \cite[Chap.~22]{AS} and \cite[Chap.~18]{NIST}.
We also note that the basic properties of the Euler $\Gamma$ function were 
used to convert the leading constant in \eqref{PS-decomposition} into that 
in \eqref{ct}. Our choice of $C_{\rho}$ is tightly connected with \eqref{identity}.
For later purposes, we let 
\begin{equation}
\label{B_rho}
B_{\rho} \define C_\rho \, \frac{\Gamma\(\frac{1+ \rho}{2}\) \, \Gamma\(1-\frac{\rho}{2}\)}{\Gamma\(\frac{1}{2}\)} < 0
\end{equation}
for $0<\rho\le 1$.

In this article we give a direct proof of the existence of a negative eigenvalue 
for the operator \eqref{operator}, then prove recursion formulas for power sums
for its eigenvalues when $\rho=1$. These power sums provide a means of 
approximating this unique negative eigenvalue. 
This problem has arisen in recent work by one of 
us \cite{SAITO-LAPEIG-ACHA} who developed the theory and applications of 
an integral operator commuting with the Laplacian defined on a general domain 
$\Omega \subset \Rf^d$, $d \geq 1$, satisfying rather interesting non-local
boundary condition.  In particular, for $d=1$, as Section~\ref{sec:nonlocal-bdry}
reviews this case in detail, the integral operator $\K_{1,1/2}$ defined 
in \eqref{operator} was shown to commute with the second order differential
operator $-\difh{2}{}{x}$ with non-local boundary condition.
In this article, we focus on the analysis of the spectra of
$\K_{\rho,a}$ for $0 < \rho \leq 1$ despite the fact that $\K_{\rho,a}$ with 
$\rho \neq 1$ does not commute with such a simple 2nd order differential
operator and that \eqref{PS-decomposition} is also valid for $-1 < \rho < 0$ 
(see Remark~\ref{kac-1d}).
The problem is certainly classical, but the results are new. We also show that techniques 
for the continuous case can be borrowed to provide new proofs for the discrete setting of distance 
matrices described in \cite{Bobo1,Bobo2}.

We let $L^2[-a,a]$ be the space of square integrable functions on 
the interval $[-a,a]$. We are interested in the following eigenvalue problem
\begin{equation} \label{evalue}
\K_{\rho, a} f(x) =\mu \, f(x) .
\end{equation}
That $\K_{\rho, a}$ has a discrete spectrum $\mu_0 \le \mu_1 \le \ldots$ is clear from the symmetry of the kernel and a simple compactness argument; viz. by the Cauchy-Schwarz inequality
\begin{equation} \label{compact}
|\K_{\rho, a} f(x)| \le |C_\rho| \, \sup_{x \in[-a,a]} \, \left(\displaystyle{\int_{-a}^{a} |x-y|^{2 \rho} \dd{y}}\right)^{1/2} \, \|f\|_2
\end{equation}
where $\|f\|_2 = \left(\int_{-a}^a f^2(x) \dd{x}\right)^{1/2}$. We are specifically interested in closed form formulas for $\sum_{n=0}^{\infty} \mu_n^p$, $p \in \N$. These are sometimes called \emph{Rayleigh functions} corresponding to the eigenvalue problem \eqref{evalue}. It is well-known that 
\begin{equation} \label{somme}
\sum_{n=0}^{\infty} \mu_n^p= \int_{-a}^{a} K_p(x,x) \dd{x}
\end{equation}
where $K_p(x,y)$ denotes the $p-$th iterated integral of $K(x,y) \define C_\rho |x-y|^{\rho}$ defined recursively
by $K_1(x,y)=K(x,y)$, and 
\begin{equation*}
K_{p+1}(x,y)= \int_{-a}^{a} K_{p}(x,z) \, K(z,y) \dd{z}, \quad p=1, 2, \ldots. 
\end{equation*}
The first couple of terms of \eqref{somme} can be directly inferred from the iteration process. For instance 

\begin{equation}
\sum_{n=0}^{\infty} \mu_n= 0
\end{equation}
and, 
\begin{equation} \label{sum}
\sum_{n=0}^{\infty} \mu_n^2= \int_{-a}^{a} K_2(x,x) \dd{x} = \left(C_\rho\right)^2 \frac{(2a)^{2(1+\rho)}}{(1+2 \rho) \, (1+\rho)}. 
\end{equation}
This paper develops explicit recursions formulas for these power sums in the limit when $\rho=1$. 

Our recursion formulas emulate those developed by various authors for Rayleigh functions, or power sums, involving roots of various transcendental equations. It was Euler who first found the first few closed expressions for what later came to be known as the Rayleigh function \cite{EULER1781} (see also \cite[Sec.~15.5]{WATSON-BESSEL}, \cite{DUTKA}):
\begin{equation}\label{sigma}
\sigma_{2\ell}(\nu) \define \sum_{n=1}^{\infty} \frac{1}{j_{\nu,n}^{2\ell}},
\quad \ell = 1, 2, \ldots,
\end{equation}
where $j_{\nu,n}$ denotes the $n$-th positive root of $z^{-\nu} J_{\nu}(z)$, and $J_{\nu}(z)$ is the Bessel function of the first kind of order $\nu$ \cite[Chap.~9]{AS}, \cite[Chap.~10]{NIST}. 
Euler's method was further developed by Lord Rayleigh \cite{RAYLEIGH}
and Carlitz \cite{CARLITZ1}. 
Both Euler and Rayleigh analyzed eigenvalues of oscillations of physical
systems (a hanging chain for Euler and a circular membrane for Rayleigh), 
which aroused their interest in computing zeros of the Bessel functions.
By exploiting a differential equation of Riccati-type satisfied by the function $z^{-\nu} J_{\nu}(z)$, Kishore \cite{KISHORE1, KISHORE2, KISHORE3} developed recursion formulas for 
$\sigma_{2\ell}(\nu)$, starting with the known expression, due to Euler and Rayleigh
\begin{eqnarray}
\label{needed}
\sigma_{2}(\nu) &=& 
\sum_{n=1}^{\infty} \frac{1}{j_{\nu,n}^{2}} = \frac{1}{4(\nu+1)} \notag \\
\sigma_{4}(\nu) &=& 
\sum_{n=1}^{\infty} \frac{1}{j_{\nu,n}^{4}} = \frac{1}{16(\nu+1)^2 (\nu+2)}. 
\end{eqnarray}

In his famous book \cite{RAYLEIGH-SOUND}, Lord Rayleigh was further led, in the context of treating the transverse vibrations of a clamped beam, to finding summation formulas for the reciprocal 4th and 8th powers of the positive roots of the equation
\begin{equation} \label{rayleigh}
\cos x \, \cosh x \pm 1 = 0.
\end{equation}
If these roots are denoted $\{m_k\}_{k=1}^{\infty}$, Lord Rayleigh found (see p.~279 of \cite{RAYLEIGH-SOUND}): 
\begin{eqnarray*}
\sum_{k=1}^{\infty} m_k^{-4} &=& \frac{1}{12} \\
\sum_{k=1}^{\infty} m_k^{-8} &=& \frac{33}{5040}.
\end{eqnarray*}

The early history of the techniques of proving these power sum formulas can be found in Watson's book \cite[Sec.~15.5]{WATSON-BESSEL} 
as well as \cite{AYOUB, DUNHAM-EULER, VARADARAJAN-BAMS}.
The more recent articles  \cite{GUPTA-MULDOON, ISMAIL-MULDOON, KERIMOV} offer 
modern views, survey recent results, and apply the techniques to various 
transcendental functions. 

Properly speaking, the technique of resolution of many of these problems goes back to Euler and his famous resolution of the ``Basel'' problem, named after the native Swiss city of Euler and the Bernoulli brothers. Euler successfully solved the problem first posed by Pietro Mengoli in 1644 \cite{DUNHAM-EULER, VARADARAJAN-BAMS} and found a closed form for the expression $\sum_{n=1}^{\infty} \frac{1}{n^2}$. It is now folklore that the sum is $\frac{\pi^2}{6}$. Heuristically, Euler's argument of 1740 \cite{EULER1740} (see also \cite{DUNHAM-EULER, VARADARAJAN-BAMS}) amounted to writing $\frac{\sin x}{x}$ in two different ways: as a Maclaurin series and as the infinite product
\[\prod_{n=1}^{\infty} \left(1-\frac{x^2}{n^2 \pi^2}\right),\]
since the roots of the transcendental equation $\sin x/x=0$ are given by $x=\pm n \pi$, for $n=1,2,\ldots$ 
Expanding the product, and equating the coefficients of $x^2$ gives the above formula. For rigorous justifications of these formulas one should consult \cite[Chap.~1]{KNOPP-COMPLEX-II}. Euler's technique is exactly what  Rayleigh employed in the case of equation \eqref{rayleigh}. Many nice examples illustrating this technique appear in the excellent paper of Speigel \cite{SPEIGEL} where generalizations of Newton's known formulas for the symmetric sums of the roots of a polynomial can be found (see also the comments in \cite{ISMAIL-MULDOON}). 

Radoux \cite{RADOUX}, Liron \cite{LIRON1, LIRON2, LIRON3}, 
and more recently Gupta-Muldoon \cite{GUPTA-MULDOON} and Ismail-Muldoon \cite{ISMAIL-MULDOON} employed similar techniques to generate various recursion formulas in the same spirit. In the case of Radoux and Liron, one finds explicit and recursive formulas for sums of even powers of reciprocals for the roots of the equation $\tan x=x$, and $\cot x=x$. 
To illustrate the case of the equation, 
$\tan x=x$, with $x_1, x_2, \ldots$ denoting the strictly positive roots of the
equation, they derived the sums of even powers of $x_k$'s, 
i.e., $\sum_{k=1}^\infty x_k^{-2 \ell}$, $\ell=1, 2, \ldots$.
For example, the cases $\ell=1, 2$ lead to
\[\sum_{k=1}^{\infty} \frac{1}{x_k^2}=\frac{1}{10}, \]
\[\sum_{k=1}^{\infty} \frac{1}{x_k^4}=\frac{1}{350}.\]
All of these are manifestations of convolution formulas relating the trace of the compact operator defined by the Green's function, and power sums of the eigenvalues as detailed in \cite{GOODWIN} and the classical book of Mikhlin \cite{MIKHLIN-INT}. The recent survey paper of Grieser \cite{GRIESER} offers a view that relates these formulas to what is known for matrices. As in \cite{GRIESER}, our work here also illustrates parallels between the continuous and discrete settings. 

Radoux \cite{RADOUX} attributes the method of finding sums of reciprocals of powers of eigenvalues of certain operators to S\`erge Nicaisse, but as detailed 
in \cite{GOODWIN, GRIESER,MIKHLIN-INT} this is truly classical.

The organization of this article is as follows.
In Section~\ref{sec:negative} we prove the existence, uniqueness, and simplicity of a negative eigenvalue for $\mathcal{K}_{\rho,a}$ directly. In Section \ref{sec:distance} we provide the means of proving the existence of this eigenvalue when dealing with distance matrices. In Sections~\ref{sec:nonlocal-bdry}-\ref{sec:recursion-formulas} we focus on the $\rho=1$ case, provide a series of standard reductions to simpler eigenvalue problems, and offer two different proofs of a recursive scheme to obtain explicit values of Rayleigh functions for \eqref{evalue} with $\rho=1$. Our main contribution in these sections are Theorems~\ref{thm:trace}, \ref{thm:highorder}, and \ref{thm:recursion}. For these sections, the proofs of the first two theorems are demonstrated directly using the properties of the eigenvalues of the non-local BVP without using the trace formulas unlike the way Goodwin proved for the regular BVPs \cite{GOODWIN}. The proof of Theorem~\ref{thm:recursion} uses the generating functions as Radoux \cite{RADOUX} and Liron \cite{LIRON1} did for different BVPs (see also Ismail and Muldoon \cite{ISMAIL-MULDOON}). Finally in Section~\ref{sec:higher}, we discuss higher dimensional considerations focusing on the centrality of the P\'olya-Szeg\H{o} expansion \eqref{PS-decomposition}.

\section{Unique Simple Negative Eigenvalue}
\label{sec:negative}

We will offer direct analytical proofs of both the existence and uniqueness of a negative eigenvalue for problem \eqref{evalue}. A probabilistic proof is offered in 
\cite{TAKASU}. We also note similar considerations in \cite{Kac57, Kac70,SPITZER1,SPITZER2}. Analytical proofs for the case of the logarithmic potential in 2-dimensions are offered in \cite{BOJD70,TROUT67}.

The fundamental eigenvalue of \eqref{evalue} is characterized by the Rayleigh-Ritz principle
\begin{equation} \label{rr}
\mu_0=\inf_{f\in L^2[-a,a]} C_\rho\, \displaystyle{\frac{\int_{-a}^{a} \int_{-a}^{a}  \, |x-y|^{\rho} f(x) f(y) \dd{x} \dd{y}}{\int_{-a}^{a} f^2(x) \dd{x}}} .
\end{equation}
\begin{prop} (Existence) \label{prop1}
The eigenvalue problem \eqref{evalue} admits at least one negative eigenvalue. 
\end{prop}

\begin{proof}
The proof of this proposition is inspired by \cite{ChoquardStubbe}.
By choosing a test function $f(x)$ appropriately, we will show that $\mu_0<0$. Let $f(x)=\chi_{[0,b]} - \chi_{[b,a]}$ where $\chi$ denotes the characteristic function of the appropriate interval and $0\le b \le a$. We will show that $b$ can be chosen to make the Dirichlet integral satisfy
\begin{equation} \label{rr2}
\psi(a,b,\rho) \define C_\rho\, \displaystyle{\int_{-a}^{a} \int_{-a}^{a}  \, |x-y|^{\rho} f(x) f(y) \dd{x} \dd{y} } < 0.
\end{equation}
This expression reduces to 
\begin{equation} \label{rr3}
\psi(a,b,\rho)=C_\rho\, \displaystyle{\int_{-a}^{a} \phi(y) f(y) \dd{y} }=C_\rho \displaystyle{\int_{0}^{b} \phi(y) f(y)  \dd{y}} - C_\rho 
\displaystyle{\int_{b}^{a} \phi(y) f(y) \dd{y} }
\end{equation}
where 
\begin{equation*}
\phi(y) \define \displaystyle{\int_{-a}^{y} (y-x)^{\rho} \, \left(\chi_{[0,b]}(x)- \chi_{[b,a]}(x) \right) \dd{x}
+\int_{y}^{a} (x-y)^{\rho} \, \left(\chi_{[0,b]}(x)- \chi_{[b,a]}(x) \right) \dd{x}}.
\end{equation*}
Simplifying further gives the expression
\begin{eqnarray*}
\phi(y)= \left\{
	\begin{array}{ll}
		-\frac{(a-y)^{\rho+1}}{\rho+1} +  2 \frac{(b-y)^{\rho+1}}{\rho+1} +\frac{(a+y)^{\rho+1}}{\rho+1}  & \mbox{for } y \in [0,b] \\
		-\frac{(a-y)^{\rho+1}}{\rho+1} -  2 \frac{(-b+y)^{\rho+1}}{\rho+1} +\frac{(a+y)^{\rho+1}}{\rho+1} & \mbox{for } y \in [b,a]
	\end{array}
\right.
\end{eqnarray*}
Performing the integrals in \eqref{rr3}, the expression in \eqref{rr2} reduces to
\begin{equation} \label{key}
\psi(a,b,\rho) =  C_\rho \displaystyle{\frac{
-(2a)^{\rho+2}+ 4 (a-b)^{\rho+2} + 2 (a+b)^{\rho+2} +2 b^{\rho+2}-2 a^{\rho+2}}
{(\rho+1) (\rho+2)}}.
\end{equation}
This is a continuous expression in $a$ and $b$. We note that since $C_\rho < 0$ for $0 < \rho \le 1$,
\begin{equation*}
\psi(a,0,\rho) = C_\rho \frac{4 a^{\rho+2} (1-2^{\rho})}{(\rho+1) (\rho+2)} > 0 \quad \text{and} \quad 
\psi(a,a,\rho) = C_\rho \frac{(2a)^{\rho+2}}{(\rho+1) (\rho+2)} < 0.
\end{equation*}
Furthermore,
\begin{equation*}
\psi(a,c a,\rho) = C_\rho \displaystyle{\frac{a^{\rho+2}}{(\rho+1) (\rho+2)}} \, \xi(c)
\end{equation*}
where $\xi(c)\define-2^{\rho+2}+ 4 (1-c)^{\rho+2}+ 2 (1+c)^{\rho+2} +2 c^{\rho+2}- 2$. This function $\xi(c)$ is monotonically increasing for $\frac{1}{2}\le c\le 1$ whereas $\psi(a, ca, \rho)$ is monotonically decreasing on the same interval, since $\xi'(c)=(\rho+2) \, \left(2 (1+c)^{\rho+1}+ 2 c^{\rho+1}- 4 (1-c)^{\rho+1} \right)>0$, viz. $1+c>c\ge 1-c$. Since $\xi(1/2)<0$, and $\xi(1)>0$, the equation $\xi(c)=0$ has a unique solution $c_0 \in (1/2, 1)$. Choose then $b$ such that 
\[c_0 a < b < a\]
to complete the proof. 
\end{proof}

\begin{rmk} \label{redux}
We note a couple of basic facts which will be useful for what follows. 

\begin{itemize}
\item[(i)] The eigenvalue problem \eqref{evalue} can be reduced to the interval $[-1,1]$ by simple rescaling.
If $\mu(\rho, a)$ and $f(x; \rho, a)$ denote an eigenvalue and the corresponding eigenfunction for $0 < \rho \leq 1$ and $a>0$, then 
\[
\mu(\rho, a) = a^{\rho+1} \mu(\rho, 1) \quad \text{and} \quad f(x; \rho, a) = f(x/a; \rho, 1), \quad x \in [-a, a].
\]
\item[(ii)] The eigenvalue problem at the origin of the investigation \cite{SAITO-LAPEIG-ACHA} was motivated by the kernel defined in \eqref{operator} with $\rho=1$ on the interval [0,1] as we shall discuss it in more detail in Section~\ref{sec:nonlocal-bdry}.  Let $\tau_\theta$ be the translation operator in $\Rf^1$ where 
$\theta \in \Rf$ defined as $\tau_\theta f(x) \define f(x-\theta)$.
Then, the integral operator, the eigenvalues, and the eigenfunctions of
this problem, denoted by $\tilde{\K}$, $\tilde{\mu}$ and $\tilde{f}$, can be 
expressed by those of $\K_{\rho, a}$, $\mu(\rho, a)$, and $f(x; \rho, a)$ as
\[
\tilde{\K} = \tau_{-\half} \K_{1,\half} \tau_{\half} ;
\quad
\tilde{\mu} = \mu\(1, \half\) = \frac{1}{4}\mu(1,1);
\quad \text{and} \quad \tilde{f}(x) = \tau_{\half} f\(x; 1, \half\).
\]
\end{itemize} 
\end{rmk}

\begin{prop} (Uniqueness) \label{prop2}
The eigenvalue problem \eqref{evalue} admits at most one negative eigenvalue. 
\end{prop}

\begin{proof}
The proof is inspired by Kac \cite{Kac57}. By virtue of Remark~\ref{redux} we will reduce the problem to the $a=1$ case. We will denote the inner product of two $L^2[-1,1]$ functions, $f, g$ by $\langle f, g \rangle := \int_{-1}^{1} f(x) g(x) \dd{x}$. 

We will prove the result by contradiction. Suppose $\mu$ and $\mu'$ are negative eigenvalues of \eqref{evalue} (not necessarily different). Let $u$ and $v$ be the corresponding eigenfunctions such that
$\langle u, v \rangle=0$. Choose $\alpha, \beta \neq 0$ such that $\langle \alpha u + \beta v, 1 \rangle=0$. Let $w= \alpha u + \beta v$. Note that $P_0^{(-\rho/2)}(x)=1$. We have 
\begin{eqnarray*}
\alpha^2 \mu+ \beta^2 \mu'&=&\int_{-1}^1 \int_{-1}^1 K(x,y) w(x) w(y) \dd{x} \dd{y} 
\notag \\
&=&
\int_{-1}^1 \int_{-1}^1 \left(K(x,y)- B_{\rho} P_0^{(-\rho/2)}(x) P_0^{(-\rho/2)}(y)\right) 
w(x) w(y) \dd{x} \dd{y}
\notag \\
&=&
B_{\rho} \, \sum_{n=1}^{\infty} \left(1- \frac{2n}{\rho}\right) \left(\int_{-1}^1 P_n^{(-\rho/2)}(x)
w(x) \dd{x}\right)^2 \quad \text{via \eqref{operator}, \eqref{PS-decomposition}, \eqref{B_rho}}
\notag \\
&\ge& 0
\end{eqnarray*}
which is a contradiction. 
\end{proof}
\begin{rmk}\label{kac-1d} 
\begin{itemize}
\item[(i)] It is instructive to compare with the case of negative powers in the kernel $K(x,y)$ (see \cite{Kac55}). In this case, there are no negative eigenvalues since for $0 < \rho < 1$, we have
\begin{eqnarray*}
\int_{-1}^1 \int_{-1}^1 \frac{1}{|x-y|^{\rho}} f(x) f(y) \dd{x} \dd{y} &=&
\frac{\Gamma\(\frac{1-\rho}{2}\) \Gamma\(1+\frac{\rho}{2}\)}{\Gamma\(\frac{1}{2}\)} \, 
\sum_{n=0}^{\infty} \left(1+ \frac{2n}{\rho}\right) \left(\int_{-1}^1 P_n^{(\rho/2)}(x)
f(x) \dd{x}\right)^2
\notag \\
&\ge& 0
\end{eqnarray*}
leading to a positive quadratic form. 
\item[(ii)] This is also the case for the 1D logarithmic potential \cite{Bekers,Reade} where all the eigenvalues are negative since the quadratic form is negative definite by virtue of the expansion 
\begin{eqnarray*}
\int_{-1}^1 \int_{-1}^1 \log|x-y| f(x) f(y) \dd{x}  \dd{y} &=& -\log 2
\left(\int_{-1}^{1} f(x) \dd{x}\right)^2-
\sum_{n=1}^{\infty} \frac{2}{n} \left(\int_{-1}^1 T_n(x)
f(x) \dd{x}\right)^2
\notag \\
&\le& 0.
\end{eqnarray*}
This follows from the well-known expansion
\begin{equation*}
\log|x-y|=-\log 2 -\sum_{n=1}^{\infty} \frac{2}{n} T_n(x) T_n(y)
\end{equation*}
where $T_n(x)$ is the Chebyshev polynomial of the first kind of order $n$. 
\end{itemize} 
\end{rmk}

\begin{prop} \label{prop3}
The operator $\mathcal{K}_{\rho,a}$ is non-singular.
\end{prop}
\begin{proof}
The proof is inspired by \cite{SPITZER1} (see also \cite{Kac70}). We will again reduce the problem to $a=1$. We will show that $\mu=0$ is not an eigenvalue. Suppose so, then $\mathcal{K}_{\rho,1} \, u=0$ for some $u\neq0$, normalized so that $\langle u, u\rangle =1$. If $\langle u, 1\rangle =0$, then as in the 
proof of Proposition~\ref{prop2},
\begin{eqnarray*}
\langle \mathcal{K}_{\rho,1} u, u \rangle&=& 
B_{\rho} \, \sum_{n=1}^{\infty} \left(1- \frac{2n}{\rho}\right) \left(\int_{-1}^1 P_n^{(-\rho/2)}(x)
u(x) \dd{x}\right)^2 
> 0.
\end{eqnarray*}
(Note that $\langle \mathcal{K}_{\rho,1} u, u \rangle=0$ means 
$\langle u, P_n^{(-\rho/2)} \rangle=0$ for $n=0,1, \ldots$. Thus $u\equiv 0$, which contradicts the fact that $u$ is an eigenfunction). 
Hence we must have $\langle u, 1\rangle \neq 0$. Let $\mu<0$ be the unique negative eigenvalue, with $\mathcal{K}_{\rho,1} v= \mu v$, and 
$\langle v, v \rangle =1$.\\ \noindent Let $\alpha, \beta \neq 0$ such that $\langle \alpha u+ \beta v, 1\rangle=0$. Again by the same argument in the proof of Proposition~\ref{prop2}, for $w=\alpha u + \beta v$, $\beta^2 \mu=\langle w,\mathcal{K}_{\rho,1} w \rangle \ge 0$ which contradicts the fact that $\mu<0$. 
Hence, $\mu=0$ is not an eigenvalue of $\K_{\rho,1}$.
\end{proof}

As a result of Proposition~\ref{prop3}, $\mathcal{A}_{\rho,a} \define \left(\mathcal{K}_{\rho,a}\right)^{-1}$ exists. Moreover, the equation 
\begin{equation}
\mathcal{A}_{\rho,a}\, v=1
\end{equation}
has a unique solution.

\begin{figure}[t]
\begin{center}
\includegraphics[width=0.85 \textwidth]{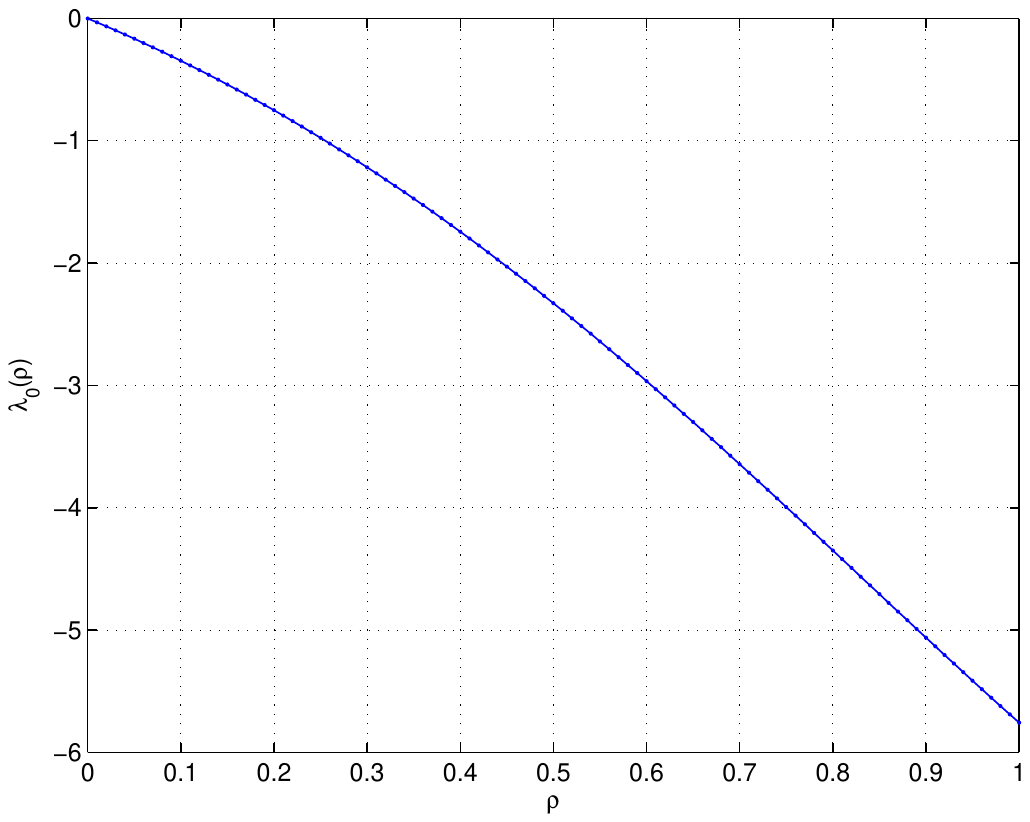}
\caption{Reciprocal negative eigenvalues, $\lambda_0(\rho) = 1/\mu(\rho,1)$, $0<\rho\le1$}\label{fig1}
\end{center}
\end{figure}

Let \[\frac{1}{R_0(\rho,a)}=\int_{-a}^{a} v(x) dx=\langle v, 1\rangle.\]
$R_0$ is the one-dimensional equivalent of the Robin constant defined in \cite{Kac70}. Its sign is tightly associated with the existence of a negative eigenvalue. This fact is exploited in \cite{BOJD70,Kac70,TAKASU} where it is demonstrated that the underlying operator has a negative eigenvalues is equivalent to $R_0<0$. We note also that $v$ has the explicit expression
\begin{equation} \label{v-expl}
v(x)=\mathcal{K}_{\rho,a} \, 1 = C_{\rho} \, \int_{-a}^{a} |x-y|^{\rho} \dd{y}
\end{equation}

Again we focus on the $a=1$ case. As before let $\mu_0,\mu_1, \ldots$ denote the eigenvalues of \eqref{evalue}, and let $u_0, u_1, \ldots$ denote the associated normalized eigenfunctions. It follows at once that 
\begin{equation} \label{robin}
\frac{1}{R_0(\rho,1)}= \sum_{n=0}^{\infty} \mu_n(\rho,1) \left(\int_{-1}^1 u_n(x) \dd{x} \right)^2
\end{equation}
(to obtain this statement, simply expand the function $v$ 
w.r.t.\ $\{ u_n \}$, then integrate). 
An explicit calculation leads to 
\[R_0(\rho,1)= \frac{(1+\rho) (2+\rho)}{2^{3+\rho} C_{\rho}}<0\]
(simply integrate \eqref{v-expl}.) 
This calculation and \eqref{robin} lead to a lower bound estimate for $\mu_0(\rho,1)$ which follows from dropping the positive terms in the series and applying the Cauchy-Schwarz inequality to $\left(\int_{-1}^1 u_0\right)^2$, namely
\begin{equation}
\mu_0(\rho,1)< \frac{2^{2+\rho} C_{\rho}}{(1+\rho) (2+\rho)} <0. 
\end{equation}

\begin{rmk}
Troutman proved in \cite{TROUT67} a similar bound for the negative eigenvalue of the logarithmic potential in terms of the transfinite diameter of the underlying domain (see also \cite{TAKASU}). In Fig.~\ref{fig1} we plot $\lambda_0(\rho) \define \frac{1}{\mu_0(\rho,1)}$ as a function of $0<\rho\le 1$. 
\end{rmk}

\begin{rmk} \label{rrr}
When $\rho=1$, the Green's function for the Dirichlet eigenvalue problem on [0,1] is given by \[G_D(x,y)=\min(x,y) - x y = \frac{1}{2} (x+y) - x y - \frac{1}{2}|x-y|\] 
which clearly indicates that our kernel is a finite-rank perturbation of the Dirichlet kernel. Indeed, 
\begin{equation*}
\tilde{\K} = \mathcal{G}_D+ \mathcal{T}_D
\end{equation*}
where $ \mathcal{G}_D$ denotes the integral operator corresponding to the Dirichlet kernel, and \[\mathcal{T}_D: L^2[0,1] \longrightarrow L^2[0,1]\] is defined by
\begin{equation*}
\mathcal{T}_D \, f(x) \define \int_0^1 \, \left(x y - \frac{1}{2} \left(x+y\right) \right) \, f(y) \dd{y}.
\end{equation*}
One can in fact calculate the eigenvalues of this perturbation. We proceed as in \cite[pp.~271--276]{invitation1}, \cite[pp.~215--216]{invitation2}. Let $u_1(x)=x$, and $u_2(x)=1$. Then, as in \cite{invitation2}, 
$\mathcal{T}_D= x_1^{\ast}(f) u_1+ x_2^{\ast}(f) u_2$
with $x_1^{\ast}(f)= \int_0^1 \left(y-\frac{1}{2}\right) f(y) \dd{y}$, 
$x_2^{\ast}(f)= \int_0^1 \(-\frac{1}{2} y \) f(y) \dd{y}$. The nonzero eigenvalues of this finite rank operator are given by the nonzero
eigenvalues of the matrix
\begin{eqnarray*} A&=& \left( \begin{array}{cc}
x_1^{\ast}(u_1) & x_1^{\ast}(u_2) \\
x_2^{\ast}(u_1) & x_2^{\ast}(u_2) 
\end{array} \right)
\notag \\
&=&
\left( \begin{array}{cc}
\int_0^1 \left(y^2-\frac{1}{2} y \right) y \dd{y} & \int_0^1 \left(y-\frac{1}{2}\right) \dd{y} \\
& \\
 - \int_0^1 \frac{1}{2} y^2 \dd{y} & - \int_0^1 \frac{1}{2} y  \dd{y}
\end{array} \right)
\notag \\
&=&
\left( \begin{array}{cc}
1/12 & 0 \\
 -1/6& -1/4 
\end{array} \right),
\end{eqnarray*}
i.e., $\lambda_1^{\ast}=-1/4$, $\lambda_{2}^{\ast}=1/12$, with corresponding eigenvectors $v_1^{\ast}(x)=-1$, $v_2^{\ast}=-(x-1/2)$. 
$\mathcal{T}_D$ is a rank 2 correction of $\tilde{\K}$ with one positive eigenvalue, and one negative eigenvalue, zero being an eigenvalue of infinite multiplicity.  Our operator is nothing but a rank 1 perturbation of a positive operator. The same arguments of 
\cite{OSELEDETS} can be applied to prove the uniqueness of a negative eigenvalue for $\tilde{\K}$.  
\end{rmk}

\vskip 1 cm

\begin{rmk} \label{rrr2}
As in Remark~\ref{rrr}, the same can be said about the Green's function  for the Neumann eigenvalue problem on [0,1], and our kernel. Since the Neumann kernel is given by
\[G_N(x,y)=-\max(x,y)+\frac{1}{2} \left(x^2+y^2\right) + \frac{1}{3}=-\frac{1}{2}|x-y|
-\frac{1}{2} \left(x+y\right)+\frac{1}{2} \left(x^2+y^2\right) + \frac{1}{3},\]
we conclude that 
\begin{equation*}
\tilde{\K} = \mathcal{G}_N+ \mathcal{T}_N
\end{equation*}
where $ \mathcal{G}_N$ denotes the integral operator corresponding to the Neumann kernel, and \[\mathcal{T}_N: L^2[0,1] \longrightarrow L^2[0,1]\] is defined by
\begin{equation*}
\mathcal{T}_N \, f(x) \define \int_0^1 \, \left(\frac{1}{2} \left(x+y\right) - \frac{1}{2} \left(x^2+y^2\right) - \frac{1}{3} \right) \, f(y) \dd{y}.
\end{equation*}
We now calculate the eigenvalues of this perturbation. Let $u_1(x)=\frac{1}{2} (x-x^2)$, and $u_2(x)=1$. Then, as in \cite{invitation2}, 
$\mathcal{T}_N= x_1^{\ast}(f) u_1+ x_2^{\ast}(f) u_2$
with $x_1^{\ast}(f)= \int_0^1 \, f(y) \dd{y}$, 
$x_2^{\ast}(f)= \int_0^1 \left(\frac{1}{2} y - \frac{1}{2} y^2 - \frac{1}{3} \right) f(y) \dd{y}$. The nonzero eigenvalues of this finite rank operator are given by the nonzero
eigenvalues of the matrix
\begin{eqnarray*} A&=& \left( \begin{array}{cc}
x_1^{\ast}(u_1) & x_1^{\ast}(u_2) \\
 x_2^{\ast}(u_1) & x_2^{\ast}(u_2) 
\end{array} \right)
\notag \\
&=&
\left( \begin{array}{cc}
\int_0^1  \frac{1}{2} \left(y-y^2\right) \dd{y} & \int_0^1 \dd{y} \\
& \\
 \int_0^1 \left(-\frac{1}{3}+\frac{1}{2} y -\frac{1}{2} y^2\right)  
 \, \frac{1}{2} \left(y-y^2\right) \dd{y} & \int_0^1 \left(-\frac{1}{3}+\frac{1}{2} y -\frac{1}{2} y^2\right)  \dd{y}
\end{array} \right)
\notag \\
&=&
\left( \begin{array}{cc}
1/12   & 1 \\
-7/360 & -1/4 
\end{array} \right),
\end{eqnarray*}
i.e., $\lambda_1^{\ast}=\frac{-5 - \sqrt{30}}{60} \approx -0.17462$, $\lambda_{2}^{\ast}=\frac{-5 + \sqrt{30}}{60} \approx 0.00795$. 
$\mathcal{T}_N$ is also a rank 2 correction of $\tilde{\K}$ with one positive eigenvalue, and one negative eigenvalue, zero being an eigenvalue of infinite multiplicity, and the same arguments of Remark~\ref{rrr} hold. 
\end{rmk}

\vskip 1 cm

\section{Distance Matrices and Matrices of Negative-Type} 
\label{sec:distance}

Much of the work of Section~\ref{sec:negative} can be emulated for the discrete case. We illustrate this in the one-dimensional case, and relegate discussion in the higher dimensional setting for an upcoming paper. 

For $-1\le x_1 < \cdots <x_{m}\le 1$, the matrix $D=\left( |x_i - x_j|^{\rho}\right)$ is called a \emph{distance matrix}. These matrices are the subject of renewed interest in recent treatments \cite{BapatKirklandNeumann,Bobo1,Bobo2,Karoui,Jiang} where $D$ is identified as a matrix of \textit{negative-type}. Such matrices have the property of exhibiting a unique simple positive eigenvalue (Multiplying by the negative coefficient $C_{\rho}$ corresponds to the results of Section~\ref{sec:negative}). A matrix $N=\left(N_{ij}\right)_{m \times m}$ is said to be of negative-type whenever the associated quadratic form 
\[\sum_{i,j=1}^m N_{ij} \xi_i \xi_j\le 0\]
for all choices $\xi_i, i=1, \ldots, m$, such that $\sum_{i=1}^m \xi_i=0$. 

We have traced the earliest works on these matrices to mid-1930s, in particular the papers of Schoenberg \cite{Schoenberg1,Schoenberg2,Schoenberg3} and Szeg\H{o} \cite{Szego} where the existence and uniqueness of this positive eigenvalue is proved directly via various techniques. Schoenberg introduced his transformation technique to allow for higher dimensional considerations while Szeg\H{o} relates the problem to Toeplitz forms. We also note that the existence and uniqueness of the positive eigenvalue follows from Perron-Frobenius theory; see in particular \cite{Gantmacher}. Much of this is reviewed and updated in the recent series of papers \cite{Bobo1,Bobo2,Jiang}. One can adapt the above propositions to the finite dimensional setting, also exploiting \eqref{PS-decomposition}. We will illustrate this connection by showing that $D$ is of negative-type, and thus by virtue of \cite{Bobo2}, it possesses only one positive simple eigenvalue. We note that for $0<\rho\le 1$
\begin{eqnarray} \label{1d-discrete}
\sum_{i, j=1}^{m} |x_i - x_j|^{\rho} \xi_i \, \xi_j &=& 
\frac{\Gamma\(\frac{1+ \rho}{2}\) \, \Gamma\(1-\frac{\rho}{2}\)}{\Gamma\(\frac{1}{2}\)} \,
\sum_{i, j}  \xi_i \, \xi_j \sum_{n=0}^{\infty} \left(1-\frac{2n}{\rho}\right) \, P_n^{(-\rho/2)} (x_i) P_n^{(-\rho/2)} (x_j)\notag \\
&=& \frac{\Gamma\(\frac{1+ \rho}{2}\) \, \Gamma\(1-\frac{\rho}{2}\)}{\Gamma\(\frac{1}{2}\)} \, \sum_{n=0}^{\infty} \left(1-\frac{2n}{\rho}\right) \, \sum_{i, j}   P_n^{(-\rho/2)} (x_i) \xi_i   P_n^{(-\rho/2)} (x_j) \xi_j \notag \\
&=& \frac{\Gamma\(\frac{1+ \rho}{2}\) \, \Gamma\(1-\frac{\rho}{2}\)}{\Gamma\(\frac{1}{2}\)} \, \sum_{n=0}^{\infty} \left(1-\frac{2 n}{\rho}\right) \, \left(\sum_{i} P_n^{(-\rho/2)} (x_i) \xi_i\right)^2. 
\end{eqnarray}
Thus, when $\sum_{i=1}^{m} \xi_i=0$, $\sum_{i, j=1}^{m} |x_i - x_j|^{\rho} \xi_i \, \xi_j\le 0$. The matrix $D$ is then 
of negative-type, and the existence and simplicity of a positive eigenvalue follows immediately from \cite{Bobo2}. One can even adapt the above arguments of Section~\ref{sec:negative} without recourse to \cite{Bobo2} or to Perron-Frobenius theory \cite[Chap.~XIII]{Gantmacher}. 

\section{A Non-local Boundary Value Problem}
\label{sec:nonlocal-bdry}
In this section, we note that problem \eqref{evalue} reduces, when $\rho=1$, with the appropriate shifts required when working on the interval $[0,1]$ as described in Remark~\ref{redux}, to the problem described in Corollary 6 from the article \cite{SAITO-LAPEIG-ACHA}, which we recall:
\begin{corollary}
\label{thm:a-line}
The eigenfunctions of the integral operator $\tilde{\K}$ with the kernel
$K(x,y) = -|x - y|/2$ for the unit interval $\Omega=(0,1)$ satisfy the following
Laplacian eigenvalue problem:
\bdm
- \phi'' = \lambda \phi, \quad x \in (0,1);
\edm
\begin{equation}
\label{eqn:nonlocal-1}
\phi(0)+\phi(1)=-\phi'(0)=\phi'(1),
\end{equation}
which can be solved explicitly as follows.
\begin{itemize}
\item $\lambda_0 \approx -5.756915$ is the smallest (and the only negative) 
eigenvalue and is the solution of the following secular equation:
\begin{equation}
\label{eqn:secular-0}
\coth \frac{\sqrt{-\lambda_0}}{2} = \frac{\sqrt{-\lambda_0}}{2},
\end{equation}
The corresponding eigenfunction is:
\begin{equation*}
\label{eqn:phi0-line}
\phi_0(x) = c_0 \cosh \sqrt{-\lambda_0} \(x-\frac{1}{2}\),
\end{equation*}
where $c_0 =\sqrt{2}\(1+\frac{\sinh\sqrt{-\lambda_0}}{\sqrt{-\lambda_0}}\)^{-1/2}$
$\approx 0.7812598$ is a normalization constant to have $\| \phi_0 \|_{L^2(\Omega)}=1$.
\item $\lambda_{2m-1}=(2m-1)^2 \pi^2$, $m=1,2,\ldots$, 
and the corresponding eigenfunction is:
\begin{equation*}
\label{eqn:phi2m-1-line}
\phi_{2m-1}(x) = \sqrt{2} \cos(2m-1)\pi x .
\end{equation*}
These are canonical cosines with odd modes.
\item $\lambda_{2m}$, $m=1,2,\ldots$, is the solution of the secular equation:
\begin{equation}
\label{eqn:secular-2m}
\cot \frac{\sqrt{\lambda_{2m}}}{2} = -\frac{\sqrt{\lambda_{2m}}}{2},
\end{equation}
and the corresponding eigenfunction is:
\begin{equation*}
\label{eqn:phi2m-line}
\phi_{2m}(x)=c_{2m} \cos \sqrt{\lambda_{2m}}\(x-\frac{1}{2}\),
\end{equation*}
where
$c_{2m}=\sqrt{2} \left\{ 1+ \frac{\sin \sqrt{\lambda_{2m}}}{\sqrt{\lambda_{2m}}} \right\}^{-1/2}$ is a normalization constant.
\end{itemize}
\end{corollary}

\begin{rmk}
We refer the reader to \cite{SAITO-LAPEIG-ACHA} for the motivation
of considering such an integral operator $\tilde{\K}$, the description of
the higher dimensional versions, and a variety of applications.
Here, however, we would like to point out our new interpretation of the above
eigenvalue problem that was not explicitly stated in \cite{SAITO-LAPEIG-ACHA}.
The above problem turns out to be equivalent to the following problem 
defined for the whole real axis and then restricting the solutions to
the unit interval $\Omega$.
\bdm
-\psi'' = \begin{cases} \lambda \psi & \text{for $x \in \Omega$}; \\
                               0     & \text{for $x \in \Rf\setminus\overline{\Omega}$},
          \end{cases}
\edm
with the continuity conditions at the boundary points: 
$\psi(0-)=\psi(0+)$, $\psi'(0-)=\psi'(0+)$,
$\psi(1-)=\psi(1+)$, $\psi'(1-)=\psi'(1+)$.
Then, $\phi(x)$ in Corollary~\ref{thm:a-line} is $\chi_\Omega (x)\psi(x)$.
\end{rmk}

\begin{rmk}
The three cases of the eigenvalues in Corollary~\ref{thm:a-line},
i.e., $\lambda_0$; $\{\lambda_{2m-1}\}$; and $\{\lambda_{2m}\}$ can also be 
derived from a single equation:
\begin{equation*}
\label{eqn:whole}
\( e^{\alpha/2} + e^{-\alpha/2} \) \cdot
\( \frac{e^{\alpha/2} + e^{-\alpha/2}}{e^{\alpha/2} - e^{-\alpha/2}} -\frac{\alpha}{2}
\) = 0,
\end{equation*}
where $\lambda=-\alpha^2$, and $\alpha \in \C$.
Searching zeros of the first factor for $\alpha \in \im \Rf$ leads to 
$\lambda_{2m-1}=(2m-1)^2\pi^2$
whereas doing so in the second factor for $\alpha \in \Rf$ leads to \eqref{eqn:secular-0} and for $\alpha \in \im \Rf$ leads to \eqref{eqn:secular-2m}.
\end{rmk}

\begin{rmk}
Both Radoux \cite{RADOUX} and Liron \cite{LIRON1} dealt with
the secular equation $\tan \beta = \beta$. They explicitly mention that
this equation came from the one-dimensional Laplacian eigenvalue problem
by setting $\lambda=-\alpha^2$, $\alpha=\im\beta$, $\beta \in \Rf$ with the 
following \emph{Robin} boundary condition:
\bdm
\phi(0) = 0, \quad \phi'(1) = \phi(1).
\edm
Note that $\phi'(0) = \phi(0)$, $\phi(1) = 0$ lead to $\tan \beta = -\beta$.

On the other hand, Radoux also dealt with the other secular equation 
$\cot \beta = \beta$ whereas Liron treated the case involving
$\cot \beta = -\beta$.  Neither of them explained
why they wanted to treat these secular equations and 
neither of them explicitly listed the corresponding boundary condition
unlike the case of $\tan \beta = \beta$.
In fact, simple computations similar to those in \cite[Sec.~4.3]{STRAUSS}
suggest that $\cot \beta = \beta$ is associated with
the Robin boundary conditions 
$(\phi'(0), \phi'(1))=(\phi(0), 0)$ 
or $(0, -\phi(1))$
and $\cot \beta = -\beta$ is associated with
$(\phi'(0), \phi'(1))=(0, \phi(1))$ 
or $(-\phi(0), 0)$.
But one also needs to consider the hyperbolic versions, i.e.,
$\coth \alpha = \alpha$, in order to fully solve the
eigenvalue problems with the Robin boundary conditions
$(\phi'(0), \phi'(1))=(0, \phi(1))$ or $(-\phi(0), 0)$.
Note that these Robin boundary conditions are all decoupled, i.e., local.
To the best of our knowledge, \cite{SAITO-LAPEIG-ACHA} is the first to explicitly
describe the unusual non-local boundary condition \eqref{eqn:nonlocal-1}.
\end{rmk}

\begin{rmk} \label{r4}
One can exploit the well-known trace formula \cite{GOODWIN,GRIESER,MIKHLIN-INT}
\begin{equation}
\label{trace}
\sum_{n=0}^\infty \frac{1}{\lambda^p_n} = \int_0^1 K_p(x,x) \dd{x},
\end{equation}
where $K_p(x,y)$ denotes the $p$th iterated kernel of $K(x,y)$, to determine the first few expressions for the Rayleigh function at hand. Indeed, one obtains at once
\[\sum_{n=0}^\infty \frac{1}{\lambda_n} = \int_0^1 K(x,x) \dd{x}=0,
\]
and
\[\sum_{n=0}^\infty \frac{1}{\lambda^2_n} = \int_0^1 K_2(x,x) \dd{x}
=\frac{1}{4} \, \int_0^1 \left(\frac{1}{3} - x + x^2 \right) \dd{x}= \frac{1}{24}.
\]
However this task becomes tedious for $p\ge 3$, and we propose to obtain these power sums without recourse to iterated kernels, but by exploiting properties of the transcendental equations of which the eigenvalues are roots. Note that this agrees with \eqref{sum} for $(\rho, a) = (1, 1/2)$. 
\end{rmk}

\section{Sum of the Reciprocals of the Eigenvalues of Corollary~\ref{thm:a-line}}
\label{sec:trace-formula}
In light of Remark~\ref{r4}, we want to show the following directly.
\begin{thm}
\label{thm:trace}
Let $\{ \lambda_{n} \}_{n=0}^\infty$ be the eigenvalues of the boundary problem in Corollary~\ref{thm:a-line}, and let $K(x,y) = -|x-y|/2$.  Then, they satisfy the following
trace formula:
\begin{equation*}
\label{eqn:trace}
\sum_{n=0}^\infty \frac{1}{\lambda_n} = \int_0^1 K(x,x) \dd{x} = 0.
\end{equation*}
\end{thm}

\begin{proof}
Let us group the eigenvalues into the three groups as indicated in Corollary~\ref{thm:a-line}:
\begin{equation}
\label{eqn:three-terms}
\sum_{n=0}^\infty \frac{1}{\lambda_n}
= \frac{1}{\lambda_0} + \sum_{m=1}^\infty \frac{1}{\lambda_{2m-1}} + 
\sum_{m=1}^\infty \frac{1}{\lambda_{2m}} .
\end{equation}
Now, the second term of the sum is:
\begin{eqnarray}
\label{eqn:term2}
\sum_{m=1}^\infty \frac{1}{\lambda_{2m-1}} &=& \sum_{m=1}^\infty \frac{1}{(2m-1)^2 \pi^2} \\ \nonumber
&=& \frac{1}{\pi^2} \sum_{m=1}^\infty \frac{1}{(2m-1)^2} \\ \nonumber
&=& \frac{1}{\pi^2} \( \sum_{m=1}^\infty \frac{1}{m^2} - \sum_{m=1}^\infty \frac{1}{(2m)^2} \) \\ \nonumber
&=& \frac{1}{\pi^2} \cdot \frac{3}{4} \sum_{m=1}^\infty \frac{1}{m^2} \\ \nonumber
&=& \frac{1}{\pi^2} \cdot \frac{3}{4} \cdot \frac{\pi^2}{6} = \frac{1}{8},
\end{eqnarray}
where we used the famous Basel problem identity $\sum_{m=1}^\infty 1/m^2 = \pi^2/6$ 
resolved by Euler \cite{EULER1740} (see also \cite{AYOUB, DUNHAM-EULER, GRIESER, VARADARAJAN-BAMS}).

As for the last term of \eqref{eqn:three-terms}, 
\begin{equation}
\label{eqn:term3}
\sum_{m=1}^\infty \frac{1}{\lambda_{2m}} = \frac{1}{4} \sum_{m=1}^\infty \frac{1}{x^2_m},
\end{equation}
where $x_m \define \sqrt{\lambda_{2m}}/2 > 0$ is the $m$th zero of the following 
transcendental equation; see \eqref{eqn:secular-2m}:
\begin{equation}
\label{eqn:tan}
\cot x = -x.
\end{equation}
To proceed to compute \eqref{eqn:term3} explicitly, let us analyze 
\eqref{eqn:tan} more deeply.
Following Radoux \cite{RADOUX}, let us first consider the following function and
its Maclaurin series expansion:
\begin{eqnarray}
\label{eqn:tan2}
(\cot x + x ) \cdot  \sin x 
&=& \cos x + x \sin x \\ \nonumber
&=& \( 1 - \frac{x^2}{2!} + \frac{x^4}{4!} - \cdots \) + x \cdot \( x - \frac{x^3}{3!} + \frac{x^5}{5!} - \cdots \)
\\ \nonumber
&=& 1 + \frac{x^2}{2} - \( \frac{1}{3!}-\frac{1}{4!} \) x^4
 + \( \frac{1}{5!}-\frac{1}{6!} \) x^6 - \cdots \\ \nonumber
&=& 1 + \frac{x^2}{2} - \frac{3}{4!} x^4 + \frac{5}{6!}x^6
- \cdots + (-1)^{k-1}\frac{2k-1}{(2k)!}x^{2k} + \cdots
\end{eqnarray}
Now, the function $\cos x + x \sin x$ can also be expanded into the following
infinite product in a manner similar to what Euler \cite{EULER1740} and 
Rayleigh \cite{RAYLEIGH} did 
(see also \cite{AYOUB, DUNHAM-EULER, GRIESER, SPEIGEL, VARADARAJAN-BAMS}):
\begin{equation}
\label{eqn:product}
\cos x + x \sin x
= \( 1 + \frac{x^2}{\alpha^2} \) \prod_{m=1}^\infty \( 1 - \frac{x^2}{x^2_m} \),
\end{equation}
where $\alpha \approx 1.19967864$ satisfies $\alpha = \coth \alpha$.

In other words, $x=\pm \im \alpha$ are the two (and only) pure imaginary roots
of $\cos x + x \sin x$.  
This can be verified as follows.
Let us seek for the pure imaginary zeros of $\cos x + x \sin x$ by
setting $x=\im y$, $y \in \Rf$.
Then, we have
\begin{eqnarray*}
\cos x + x \sin x &=& \cos(\im y) + \im y \sin(\im y) \\ \nonumber
&=& \frac{\e^{\im (\im y)} + \e^{-\im (\im y)}}{2} + \im y \frac{\e^{\im (\im y)} - \e^{-\im (\im y)}}{2\im} \\ \nonumber
&=& \frac{\e^{y}+\e^{-y}}{2} - y \frac{\e^{y}-\e^{-y}}{2}  = 0,
\end{eqnarray*}
which is equivalent to $\cosh y -y \sinh y = 0$, i.e.,
\begin{equation}
\label{eqn:coth}
y = \coth y.
\end{equation}
The justification for the product formula \eqref{eqn:product} follows considerations similar to those for example in \cite[Chap.~1]{KNOPP-COMPLEX-II}; see also \cite{ISMAIL-MULDOON}. From \eqref{eqn:product}, we have
\begin{eqnarray}
\label{eqn:product2}
\cos x + x \sin x &=& \prod_{m=1}^\infty \( 1 - \frac{x^2}{x^2_m} \)
+ \frac{x^2}{\alpha^2} \prod_{m=1}^\infty \( 1 - \frac{x^2}{x^2_m} \) \\ \nonumber
&=& 1 + \( \frac{1}{\alpha^2} - \sum_{m=1}^\infty \frac{1}{x^2_m} \) x^2
+ \cdots
\end{eqnarray}
Equating the corresponding coefficients of the $x^2$ terms of \eqref{eqn:tan2}
and \eqref{eqn:product2}, we have
\begin{equation*}
\label{eqn:x2}
\sum_{m=1}^\infty \frac{1}{x^2_m} = \frac{1}{\alpha^2} - \frac{1}{2}.
\end{equation*}
Hence, inserting this to \eqref{eqn:term3}, in turn, \eqref{eqn:three-terms}
together with \eqref{eqn:term2} gives us
\begin{eqnarray*}
\label{eqn:three-terms2}
\sum_{n=0}^\infty \frac{1}{\lambda_n} &=&
\frac{1}{\lambda_0} + \sum_{m=1}^\infty \frac{1}{\lambda_{2m-1}} + 
\sum_{m=1}^\infty \frac{1}{\lambda_{2m}} \\ \nonumber
&=& \frac{1}{\lambda_0} + \frac{1}{8} + \frac{1}{4} \sum_{m=1}^\infty \frac{1}{x^2_m} \\ \nonumber
&=& \frac{1}{\lambda_0} + \frac{1}{8} + \frac{1}{4} \( \frac{1}{\alpha^2} - \frac{1}{2} \) \\ \nonumber
&=& \frac{1}{\lambda_0} + \frac{1}{4\alpha^2} \\ \nonumber
&=& 0,
\end{eqnarray*}
since $\lambda_0 = -4 \alpha^2$, which can be verified by identifying
\eqref{eqn:secular-0} with the equation \eqref{eqn:coth}
via $\alpha=\sqrt{-\lambda_0}/2$.
\qquad\end{proof}

\section{Sums of Higher Powers of the Reciprocals of the Eigenvalues of Corollary~\ref{thm:a-line}}
\label{sec:gen-trace-formula}
Furthermore, we can establish the following identities:

\begin{thm}
\label{thm:highorder}
Let $\{ \lambda_n \}_{n=0}^\infty$ be the eigenvalues of the boundary value
problem specified in Corollary~\ref{thm:a-line}.
Let $K_p(x,y)$ be the $p$th iterated kernel of $K(x,y)=-|x-y|/2$.
Then, we have
\begin{equation}
\label{eqn:highorder}
\sum_{n=0}^\infty \frac{1}{\lambda^p_n} = \int_0^1 K_p(x,x) \dd{x} = 
\frac{1}{4^p} \( S_{2p} + \frac{(-1)^p}{\alpha^{2p}} \) +
\frac{4^p-1}{2 \cdot (2p)!} |B_{2p}|,
\end{equation}
where
\begin{equation*}
\label{eqn:s_2p}
S_{2p} \define \sum_{m=1}^\infty \frac{1}{x_m^{2p}} = 
\sum_{m=1}^\infty \( \frac{4}{\lambda_{2m}} \)^p,
\end{equation*}
and $B_{2p}$ is the Bernoulli number, which is defined via the generating
function:
\bdm
\frac{x}{\e^x-1} = \sum_{n=0}^\infty \frac{B_n}{n!}x^n .
\edm
Moreover, $S_{2p}$ satisfies the following recursion formula:
\begin{equation}
\label{eqn:Sn}
\sum_{\ell=1}^{n+1} \frac{(-1)^{n-\ell+1} \( 2\(n-\ell+1\)-1\)}{\( 2\(n-\ell+1\)\)!}
\left\{ S_{2\ell} + \frac{(-1)^\ell}{\alpha^{2\ell}} \right\}
= \frac{(-1)^n}{2(2n)!}.
\end{equation}
\end{thm}

\begin{proof}
The first equality in \eqref{eqn:highorder} connecting the sum of the
powers of the eigenvalues and the trace of the iterated kernel is the 
standard fact and its proof can be found in, e.g., \cite[Sec.~15]{MIKHLIN-INT}.
Now, to prove the second equality, we have
\begin{eqnarray*}
\sum_{n=0}^\infty \frac{1}{\lambda^p_n} &=&
\frac{1}{\lambda_0^p} + \sum_{m=1}^\infty \frac{1}{\lambda^p_{2m-1}}
+ \sum_{m=1}^\infty \frac{1}{\lambda^p_{2m}} \\
&=& \( \frac{-1}{4\alpha^2} \)^p 
 + \frac{1}{\pi^{2p}} \sum_{m=1}^\infty \frac{1}{(2m-1)^{2p}}
 + \frac{1}{4^p} \sum_{m=1}^\infty \frac{1}{x_m^{2p}} \\
&=& \frac{(-1)^p}{4^p\alpha^{2p}} + \frac{1}{\pi^{2p}}\( 1-\frac{1}{2^{2p}} \)
\sum_{m=1}^\infty \frac{1}{m^{2p}} + \frac{1}{4^p}S_{2p} \\
&=& \frac{1}{4^p}\left\{ S_{2p}+\frac{(-1)^p}{\alpha^{2p}} + \frac{4^p-1}{\pi^{2p}}
\sum_{m=1}^\infty \frac{1}{m^{2p}} \right\} \\
&=& \frac{1}{4^p}\( S_{2p}+\frac{(-1)^p}{\alpha^{2p}} \)
+ \frac{4^p-1}{2(2p)!} | B_{2p} |,
\end{eqnarray*}
where we used the following well-known formula first obtained by Euler
(see, e.g., \cite{AYOUB, DUNHAM-EULER, VARADARAJAN-BAMS} to derive the last equality:
\bdm
\sum_{m=1}^\infty \frac{1}{m^{2p}} = \frac{(2\pi)^{2p}}{2(2p)!} | B_{2p} |.
\edm

Now, to prove the recursion formula \eqref{eqn:Sn}, we follow 
Radoux \cite{RADOUX} again.  Taking the logarithm of the product formula
\eqref{eqn:product} followed by differentiation with respect to $x$,
we have
\bdm
\frac{x \cos x}{\cos x + x \sin x} = 
\frac{\frac{2x}{\alpha}}{1+\frac{x^2}{\alpha^2}} +
\sum_{m=1}^\infty \frac{\frac{-2x}{x_m^2}}{1-\frac{x^2}{x_m^2}},
\edm
which leads to
\begin{equation} \notag
\frac{1}{2}\cos x = (\cos x + x \sin x) \cdot
\left\{ \frac{1}{\alpha^2}\frac{1}{1+ \frac{x^2}{\alpha^2}}
-\sum_{m=1}^\infty \frac{1}{x_m^2}\frac{1}{1-\frac{x^2}{x_m^2}} \right\}.
\end{equation}
Expanding each term into the Maclaurin series or the geometric series, we have
\begin{equation} \notag
\frac{1}{2} \sum_{n=0}^\infty \frac{(-1)^n x^{2n}}{(2n)!}
= \( \sum_{k=0}^\infty \frac{(-1)^{k-1} (2k-1)}{(2k)!}x^{2k} \) \cdot
\( \sum_{\ell=0}^\infty \( \frac{(-1)^\ell}{\alpha^{2\ell+2}} - S_{2\ell+2} \) x^{2\ell} \).
\end{equation}
Hence, comparing the coefficients of the $x^{2n}$ term, we have:
\begin{eqnarray*}
\frac{(-1)^n}{2(2n)!} &=& \sum_{k=0}^n \frac{(-1)^{n-k-1}(2n-2k-1)}{(2n-2k)!}
\( \frac{(-1)^k}{\alpha^{2k+2}}-S_{2k+2} \) \\
&=& \sum_{\ell=1}^{n+1} \frac{(-1)^{n-\ell+1}(2(n-\ell+1)-1)}{(2(n-\ell+1))!}
\(S_{2\ell} + \frac{(-1)^\ell}{\alpha^{2\ell}} \) \quad \text{via setting $\ell=k+1$,}
\end{eqnarray*}
which is \eqref{eqn:Sn}.
\qquad\end{proof}

Let $A_p \define \sum_{n=0}^\infty \frac{1}{\lambda^p_n}$.
Here are the first few sums:
\begin{equation} \notag
A_1 = 0; \quad A_2=\frac{1}{24}; \quad A_3= -\frac{1}{240}, \ldots
\end{equation} 

\section{The Generating Function and Obtaining Recursive Formulas All at Once}
\label{sec:recursion-formulas}
In this section, we show how to obtain the recursion formulas for the $A_p$'s at once and without recourse to the knowledge of Bernoulli numbers. The main result is the following theorem.
\begin{thm}
\label{thm:recursion}
Let $\{ \lambda_{n} \}_{n=0}^\infty$ be the eigenvalues of the boundary problem in Corollary~\ref{thm:a-line}, and let
$K_p(x,y)$ be the $p$th iterated kernel of $K(x,y) = -|x-y|/2$.  Then, 
\[A_{p}= \sum_{n=0}^\infty \frac{1}{\lambda_n^{p}} = \int_0^1 K_p(x,x) \dd{x}\] 
satisfies the recursion formula:
\begin{equation*}
\label{eqn:recursion}
4 A_{p+1}+ \sum_{k=1}^{p-1} (-1)^k \left(
\frac{2}{(2k)!} 
- \frac{1}{(2k-1)!}\right) \, A_{p-k+1} = \dfrac{(-1)^{p+1}  p}{(2p+1)!}, \quad
p=1, 2, \dots,
\end{equation*}
with $A_1=0$.
\end{thm}

\begin{proof}
From the statement of Corollary~\ref{thm:a-line} and \eqref{eqn:product}, it is clear that
\begin{equation*} 
\label{gen.1}
\left( \cos x + x \sin x \right) \cdot \cos x  = \left(1+ \frac{x^2}{\alpha^2} \right) \, \displaystyle\prod_{m=1}^{\infty} \left(1-\frac{x^2}{x_m^2}\right)
\end{equation*} 
where $\lambda_0=-4 \alpha^2$ as defined above and where we set $x_k=\sqrt{\lambda_k}/2$, for $k=1, 2, \ldots$. One can again
justify this product formula as in Knopp \cite[Chap.~1]{KNOPP-COMPLEX-II} or any standard Complex Analysis textbook which treats the Weierstrass Factor Theorem. 

In terms of the eigenvalues one has, after some trigonometric substitutions, 
\begin{equation} 
\label{gen.2}
\frac{1 + \cos x}{2} + \frac{x}{4} \, \sin x  = \left(1-  \frac{x^2}{\lambda_0} \right) \, \displaystyle\prod_{m=1}^{\infty} \left(1-\frac{x^2}{\lambda_m}\right).
\end{equation} 

Expanding the LHS into a Maclaurin series and equating lead to

\begin{equation*} 
\label{gen.3}
1+ \sum_{k=1}^{\infty} (-1)^k
\left (\frac{1}{2 (2k)!}-\frac{1}{4 (2 k-1)!} \right) \, x^{2k}
= 1- \left(\sum_{k=0}^{\infty} \frac{1}{\lambda_k} \right) x^2
+ \left(\sum_{j,k=0}^{\infty} \, \frac{1}{\lambda_k \lambda_j} \right) x^4-\ldots. 
\end{equation*} 
With $\alpha_k \define (-1)^k \left(\frac{1}{2 (2k)!}-\frac{1}{4 (2 k-1)!}\right)$ denoting the coefficients of the Maclaurin expansion, one can recourse to Speigel's formulas \cite{SPEIGEL, GOODWIN} 
\begin{eqnarray*} 
\label{gen.4}
\sum_{k=0}^{\infty} \frac{1}{\lambda_k} &=& - \alpha_1 \notag \\
\sum_{k=0}^{\infty} \frac{1}{\lambda_k^2} &=& \alpha_1^2- 2 \alpha_2 \notag \\
\sum_{k=0}^{\infty} \frac{1}{\lambda_k^3} &=& 3 \alpha_1 \alpha_2 - 3 \alpha_3 - \alpha_1^3
\notag \\
\sum_{k=0}^{\infty} \frac{1}{\lambda_k^4} &=& \alpha_1^4 - 4 \alpha_1^2 \alpha_2 + 2 \alpha_2^2 + 4 \alpha_1 \alpha_3 - 4 \alpha_4 \notag \\
\end{eqnarray*} 
to obtain, as above, 
\begin{eqnarray*} 
\label{gen.5}
A_1=\sum_{k=0}^{\infty} \frac{1}{\lambda_k} &=& 0\notag \\
A_2=\sum_{k=0}^{\infty} \frac{1}{\lambda_k^2} &=& \frac{1}{24}\notag \\
A_3=\sum_{k=0}^{\infty} \frac{1}{\lambda_k^3} &=& -\frac{1}{240}\notag \\
A_4=\sum_{k=0}^{\infty} \frac{1}{\lambda_k^4} &=& \frac{41}{40320}. 
\end{eqnarray*} 

One can generate a recursion formula for the $A_p$ sequence employing what Ismail and Muldoon \cite{ISMAIL-MULDOON} call, properly, the ``Euler-Rayleigh'' technique. 
The logarithmic derivative of the entire function 
$f(z)=\frac{1+\cos z}{2} + \frac{z}{4} \, \sin z $ 
appearing in \eqref{gen.2} gives,

\begin{equation*} 
\label{gen.6}
\dfrac{ -\dfrac{\sin z}{4} + \dfrac{z \cos z}{4} }
{\dfrac{1+\cos z}{2} + \dfrac{z \sin z}{4}}
= 
- \dfrac{2 z}{\lambda_0 - z^2} - 2 \sum_{k=1}^{\infty} \dfrac{z}{\lambda_k - z^2}
\end{equation*} 
Or, substituting $\lambda_0 = - 4 \alpha^2$ and $\lambda_k=4 x_k^2$, and after some manipulation, 
\begin{equation} 
\label{gen.7}
\dfrac{ - \dfrac{\sin 2z}{4} + \dfrac{z \cos 2z}{2} }
{ \dfrac{1+\cos 2z}{2} + \dfrac{z \sin 2z}{2}}
= 
\dfrac{z}{\alpha^2 + z^2} - \sum_{k=1}^{\infty} \dfrac{z}{x_k^2 - z^2} =:-z G(z).
\end{equation}
The function
\[G(t)=-\dfrac{1}{\alpha^2+ t^2} + \sum_{k=1}^{\infty} \dfrac{1}{x_k^2-t^2} \]
is known as the generating function of $A_p$. That is, one can obtain the needed recursion formula for this sequence from consideration of this function. To simplify notation, we let
$M_{\ell} \define 4^{\ell+1} A_{\ell+1}$. It is then clear that
\begin{equation*}
\label{gen.8}
M_{\ell-1}=\frac{(-1)^{\ell}}{\alpha^{2 \ell}}+ \sum_{m=1}^{\infty} \frac{1}{x_m^{2 \ell}}.
\end{equation*}
Moreover, a straightforward calculation leads to
\begin{equation*}
\label{gen.9}
\sum_{\ell=0}^{\infty} M_{\ell} t^{2 \ell}=G(t).
\end{equation*}
By \eqref{gen.7}, one then obtains
\begin{equation*}
\label{gen.10}
\dfrac{\sin 2 t}{4 t}- \dfrac{\cos 2 t}{2} =
\left( \dfrac{1 + \cos 2t}{2} + \dfrac{t}{2} \sin 2t \right)
\, \left(
\sum_{\ell=0}^{\infty} M_{\ell} t^{2 \ell}
\right).
\end{equation*}
Expanding into power series leads to
\begin{equation*}
\label{gen.11}
\sum_{n=1}^{\infty} (-1)^{n+1} \dfrac{4^n n}{(2n+1)!} \, t^{2n}
=\left(1+
\sum_{k=1}^{\infty} (-1)^k \left(
\frac{2^{2k-1}}{(2k)!} 
- \frac{2^{2k-2}}{(2k-1)!}\right) t^{2k}  \right)
\, \left(
\sum_{\ell=0}^{\infty} M_{\ell} t^{2 \ell}
\right).
\end{equation*}
From which one obtains $M_0=0$, and
\begin{equation*}
\label{gen.12}
\sum_{k+\ell=p} (-1)^k \left(
\frac{2^{2k-1}}{(2k)!} 
- \frac{2^{2k-2}}{(2k-1)!}\right) \, M_{\ell} = \dfrac{(-1)^{p+1} 4^p p}{(2p+1)!}.
\end{equation*}
In terms of the $A_p$'s one has,
$A_1=0$, as before,
\begin{equation*}
\label{gen.13}
4 A_{p+1}+ \sum_{k=1}^{p-1} (-1)^k 4^{-k+1} \left(
\frac{2^{2k-1}}{(2k)!} 
- \frac{2^{2k-2}}{(2k-1)!}\right) \, A_{p-k+1} = \dfrac{(-1)^{p+1}  p}{(2p+1)!},
\end{equation*}
which is the same as the desired statement of the theorem. 
\qquad\end{proof}

\begin{rmk}
We note that the recursion generates the following values$A_2=1/24$, $A_3=-1/240$, $A_4=41/40320$, 
$A_5=-107/725760$, etc., corresponding to what we obtained differently in 
Section~\ref{sec:gen-trace-formula}.
\end{rmk}
\begin{rmk}
As in \cite{ISMAIL-MULDOON}, one can exploit the formulas generated for the $A_p$'s to obtain
\begin{equation}
\label{gen.14}
-|A_{2m-1}|^{-1/(2m-1)}<\lambda_0<- A_{2m}^{-1/(2m)}
\end{equation}
and
\begin{equation}
\label{gen.15}
A_{2m}/A_{2m+1}<\lambda_0<A_{2m-1}/A_{2m}.
\end{equation}
for $m=1, 2, 3, \ldots$. These inequalities provide strict improvable bounds for the unique negative root of the transcendental equation \eqref{eqn:secular-0} and another way of obtaining it.  
\end{rmk}

\section{Higher Dimensional Considerations}
\label{sec:higher}
One of the motivations that led to the non-local BVP considered in \cite{SAITO-LAPEIG-ACHA} is that one is able to read the spectral data (eigenvalues, eigenfunctions) by discretizing then computing integrals involving the kernel $K({\bm x},{\bm y})$ over a domain $\Omega\subset \mathbb{R}^d$ without imposing conditions on $\partial \Omega$. For the two-dimensional case, $K({\bm x},{\bm y})$ takes the form of a logarithmic kernel
\begin{equation}
\label{ker}
K({\bm x},{\bm y})=-\frac{1}{2 \pi} \log \|{\bm x}-{\bm y}\|.\end{equation}
Troutman \cite{TROUT67} gave an analytical proof  for the existence of at most one negative eigenvalue and gave an upper bound estimate for it in terms of the area and transfinite diameter of $\Omega$. (The transfinite diameter is a measure of the compactness of a domain; see \cite{TROUT67} for the definition.) In \cite{Kac70}, Kac offers a probabilistic proof of this fact (see also 
\cite{BOJD70, SPITZER1, SPITZER2} and the generalization in \cite{TAKASU}). Related works are also offered in \cite{Bekers, OSELEDETS, Reade}. 

With $A_p$ denoting the power sum in \eqref{trace} and the iterated integrals computed numerically, \eqref{gen.14} and \eqref{gen.15} provide a practical and improvable means of computing this negative eigenvalue for a specific domain. When the transfinite diameter of $\Omega$ is less than or equal to one, this negative eigenvalue disappears. This is the case of the unit disk. In \cite{SAITO-LAPEIG-ACHA}, it was found that the eigenvalues of the nonlocal BVP associated with the kernel \eqref{ker} are of two types,
$j_{0,n}^2$, with multiplicity 3, and $j_{m-1,n}^2$ with multiplicity 2,
for $m=2,3, \ldots$, and $n=1,2, \ldots$. Based on the values of the Rayleigh function $\sigma_{2p}(\nu)$ defined in \eqref{sigma}, one can generate for the first few power sums. 
While $\sum_{k=1}^{\infty} \, 1/\lambda_k$ is easily seen to diverge, we have
\begin{equation*}
\sum_{k=1}^{\infty} \, \frac{1}{\lambda_k^2} = 3 \sigma_4(0)+2 \sum_{\nu=1}^{\infty} \sigma_4(\nu) = \frac{3}{32} + \frac{1}{8} \left(\frac{\pi^2}{6}-\frac{3}{2}\right).
\end{equation*}
Similarly 
\begin{equation*}
\sum_{k=1}^{\infty} \, \frac{1}{\lambda_k^p} = 3 \sigma_{2p}(0)+2 \sum_{\nu=1}^{\infty} \sigma_{2p}(\nu)
\end{equation*}
can be carried out explicitly for $p=3, 4, \ldots$, but there may not be an obvious recursion scheme.

When the kernel takes the form 
\begin{equation}
\label{ker2}
K({\bm x},{\bm y})=\|{\bm x}-{\bm y}\|^{\rho}
\end{equation}
for $0<\rho\le 1$ on $\Omega \subset \mathbb{R}^2$ one can prove the existence of a negative eigenvalue based on the formula
\eqref{PS-decomposition} independently of classical proofs based on the Schoenberg transformation \cite{Schoenberg1,Schoenberg2,Schoenberg3}. Renewed interest focuses on the discrete case, namely that of nature of the spectrum of distance matrices  \cite{Bobo1,Bobo2} (see also \cite{Karoui,Jiang} where the density of states is treated and its limiting distribution when the size of the matrix goes to $\infty$ is determined).  We describe the procedure for $\Omega \subset \{ {\bm x} \in \mathbb{R}^2, \text{ such that } \|{\bm x}\|<1\}$. We first note the identity
\begin{equation} \label{basic2d}
\|{\bm x}\|^{\rho}=\frac{1}{C_{\rho}} \, \int_{-\pi}^{\pi} \left|{\bm x} \cdot {\bm \xi}\right|^{\rho} \dd{\sigma({\bm \xi})}
\end{equation}
where $\dd{\sigma({\bm \xi})}=\dd{\theta}$ denotes the element of arclength, and this time
\begin{equation}
C_{\rho} \define \int_{-\pi}^{\pi} |\cos \theta|^{\rho} \dd{\theta} = \frac{2 \Gamma\(\frac{1}{2}\)\, \Gamma\(\frac{1+\rho}{2}\)}{\Gamma\(1+\frac{\rho}{2}\)} > 0.
\end{equation}
Combining \eqref{basic2d} and \eqref{PS-decomposition}, one obtains
\begin{eqnarray}
\|{\bm x}- {\bm y}\|^{\rho}&=& \frac{\Gamma\(1+\frac{\rho}{2}\)}{2 \Gamma\(\frac{1}{2}\)\, \Gamma\(\frac{1+\rho}{2}\)} \, \int_{-\pi}^{\pi} |{\bm x} \cdot {\bm \xi}- {\bm y} \cdot {\bm \xi}|^{\rho} \dd{\sigma({\bm \xi})} 
\notag \\
&=& \frac{\Gamma\(1+\frac{\rho}{2}\) \, \Gamma\(1-\frac{\rho}{2}\)}{2 \pi} \, \sum_{n=0}^{\infty} \left(1-\frac{2 n}{\rho} \right)\int_{-\pi}^{\pi} P_n^{\(-\frac{\rho}{2}\)} \({\bm x} \cdot {\bm \xi}\) 
P_n^{\(-\frac{\rho}{2}\)} \({\bm y} \cdot {\bm \xi}\)  \dd{\sigma({\bm \xi})} 
\end{eqnarray}
It then becomes transparent how to proceed in the case of the quadratic form with kernel \eqref{ker2}, viz., 
\begin{equation}
\int_{\Omega} \int_{\Omega} K({\bm x},{\bm y}) f({\bm x}) f({\bm y}) \dd{{\bm x}} \dd{{\bm y}}=  \frac{\Gamma\(1+\frac{\rho}{2}\) \, \Gamma\(1-\frac{\rho}{2}\)}{2 \pi} \, \sum_{n=0}^{\infty} \left(1-\frac{2n }{\rho} \right)\int_{-\pi}^{\pi}  \left(
\int_{\Omega} f({\bm x}) P_n^{\(-\frac{\rho}{2}\)} \({\bm x} \cdot {\bm \xi}\) \dd{{\bm x}}
\right)^2 
\, \dd{\sigma({\bm \xi})}. 
\end{equation}
When $\int_{\Omega} f({\bm x}) \dd{\bm x}=0$, the quadratic form is such that 
\begin{equation*}
\int_{\Omega} \int_{\Omega} K({\bm x},{\bm y}) f({\bm x}) f({\bm y}) \dd{{\bm x}} \dd{{\bm y}}\le 0.
\end{equation*}
One can even introduce the notion of a kernel of negative-type. In the discrete case, the 2-dimensional version of \eqref{1d-discrete}, when $\{{\bm x}_i\}$ are confined to the unit disk, takes the form
\begin{equation}
\sum_{i, j} \|{\bm x}_i -{\bm x}_j\|^{\rho} t_i t_j =  \frac{\Gamma\(1+\frac{\rho}{2}\) \, \Gamma\(1-\frac{\rho}{2}\)}{2 \pi} \, \sum_{n=0}^{\infty} \left(1-\frac{2 n}{\rho} \right)\int_{-\pi}^{\pi} \, \left(\sum_i t_i  P_n^{\(-\frac{\rho}{2}\)} \({\bm x}_i \cdot {\bm \xi}\)
\right)^2 
\, \dd{\sigma({\bm \xi})}. 
\end{equation}
When $\sum t_i=0$, the quadratic form is such that 
\begin{equation*}
\sum_{i, j} \|{\bm x}_i -{\bm x}_j\|^{\rho} t_i t_j\le 0.
\end{equation*}
Thus the matrix $\left(\|{\bm x}_i -{\bm x}_j\|^{\rho}\right)$ is also of negative-type, and the results of \cite{Bobo2} can be used to complete the proof for the existence, uniqueness, and simplicity of a positive eigenvalue.

\section*{Acknowledgment}
L.H. would like to thank the UC Davis Department of Mathematics for hospitality
and support while doing some of this work.
N.S.'s research was partially supported by the ONR grants
N00014-07-1-0166, N00014-09-1-0041, N00014-09-1-0318,
and N00014-12-1-0177.


\end{document}